\documentclass[12pt,reqno,sumlimits]{amsart}

\usepackage{amssymb,amscd,amsmath,epsfig}
\usepackage{
  array,
  booktabs,
  dcolumn,
  rotating,
  shortvrb,
  tabularx,
  units,
  url,
}
\theoremstyle{definition}

\newcommand{\R}{{\mathbb R}}
\newcommand{\G}{{\mathbb G}}
\newcommand{\Q}{{\mathbb Q}}
\newcommand{\C}{{\mathbb C}}

\newcommand{\Z}{{\mathbb Z}}

\newcommand{\N}{{\mathbb N}}

\newcommand{\calA}{{\mathcal A}}

\newcommand{\calC}{{\mathcal C}}
\newcommand{\calD}{{\mathcal D}}
\newcommand{\calE}{{\mathcal E}}

\newcommand{\calG}{{\mathcal G}}
\newcommand{\calH}{{\mathcal H}}

\newcommand{\calM}{{\mathcal M}}
\newcommand{\calO}{{\mathcal O}}
\newcommand{\calP}{{\mathcal P}}
\newcommand{\calR}{{\mathcal R}}

\renewcommand{\to}{\longrightarrow}

\newcommand{\ev}{\operatorname{ev}}

\newcommand{\Tr}{\operatorname{Tr}}
\newcommand{\tr}{\operatorname{tr}}

\newsavebox{\savepar}

\numberwithin{equation}{section}

\newcounter{labelflag} 
\newcommand{\labelon}{\setcounter{labelflag}{1}}
\newcommand{\Label}[1]{
                       \ifnum\thelabelflag=1
                          \ifmmode
                             \makebox[0in][l]{\qquad\fbox{\rm#1}}
                          \else
                             \marginpar{\vspace{0.7\baselineskip}
                                        \hspace{-1.1\textwidth}
                                        \fbox{\rm#1}}
                          \fi
                       \fi
                       \label{#1}
                      }
\labelon

\usepackage[all]{xy}

\newtheorem{theorem}{Theorem}[section]
\newtheorem{lemma}[theorem]{Lemma}

\newtheorem{proposition}[theorem]{Proposition}

\theoremstyle{definition}
\newtheorem{definition}[theorem]{Definition}

\newtheorem{conjecture}[theorem]{Conjecture}

\theoremstyle{remark}
\newtheorem{remark}[theorem]{Remark}

\newcommand{\BbC}{{\mathbb C}}

\newcommand{\pdo}{\Psi{\rm DO}}
\newcommand{\pdoo}{\Psi_{\leq 0}}
\newcommand{\pdos}{\Psi_{0}^*}
\newcommand{\calg}{{\mathfrak g}}
\newcommand{\frakg}{{\mathfrak g}}

\newcommand{\dvol}{{\rm dvol}}

\newcommand{\eend}{{\rm End}}
\newcommand{\ddd}{\partial \! \! \! /}
\newcommand{\indd}{{\rm IND}}

\newcommand{\cw}{\rm cw}
\newcommand{\mapsnm}{{\rm Maps}(N,M)}
\newcommand{\mapsfnm}{{\rm Maps}_f(N,M)}
\newcommand{\resw}{{\rm res}^{\rm w}}
\newcommand{\ckw}{c_k^{\rm w}}
\newcommand{\cklo}{c_k^{\rm lo}}
\newcommand{\cswk}{CS^{\rm w}_k}
\newcommand{\dg}{\dot\gamma}
\newcommand{{\diff}}{\rm Diff}
\newcommand{\adg}{{\rm Ad}_G}
\newcommand{\csw}{CS^{\rm w}}

\begin{document}


\normalsize

\title{Chern--Weil theory for certain infinite-dimensional Lie groups}
\author[S. Rosenberg]{Steven Rosenberg${}^1$}\footnote{Department of Mathematics and Statistics, Boston University, Boston, MA 02215, USA. Email: sr@math.bu.edu.}


\maketitle

\centerline{To Joe Wolf, with great appreciation}

\begin{abstract}  Chern--Weil and Chern--Simons theory extend to certain infinite-rank bundles that
appear in mathematical physics.  We discuss what is known of the invariant theory of the corresponding infinite-dimensional Lie groups.  We use these techniques to detect cohomology classes for spaces of maps between manifolds and for diffeomorphism groups of manifolds.
\end{abstract}

\medskip


\section{{\bf Introduction}}  

The purpose of this survey is to emphasize algebraic aspects of Chern--Weil 
and Chern--Simons
theory in
infinite dimensions.  The main
open questions concern the classification or even existence of nontrivial Ad-invariant functions
on the Lie algebra of some infinite-dimensional Lie groups that naturally appear in differential geometry and mathematical physics.   These groups include gauge groups
of vector bundles, diffeomorphism groups of manifolds, groups of bounded invertible pseudodifferential operators, and semidirect products of these groups.  They typically
appear as  the structure groups of manifolds of maps between manifolds (e.g., in string theory) and in the setup of the Atiyah--Singer families index theorem.
In
effect, this  article is a plea by a geometer for help from the experts in Lie groups.  
\medskip

As reviewed in Section 2, 
Chern--Weil theory is a well-established procedure to pass from an Ad-invariant polynomial 
$p$ on the Lie algebra of a finite-dimensional Lie group $G$ defined over a field $k$ to an element  $c_p\in H^*(BG,k)$, the cohomology of the classifying space $BG$.  
For the classical compact connected groups over $\C$, this correspondence is an isomorphism.
For a $G$-bundle $P\to B$ classified by a map $f:B\to BG$, the 
class $c_p(P) = f^*c_p\in H^*(B,\C)$ is by definition the characteristic class
of $P$ associated to $p$.  For $G = U(n),  SO(n)$, 
these are the Chern classes and Pontrjagin classes, respectively.  They are used extensively in differential
geometry, algebraic geometry and differential topology.

Characteristic classes  are obstructions to bundle triviality, i.e., they vanish on trivial bundles.   
When a Chern class of a bundle vanishes, the precise obstruction information for that class
is unavailable, but
 there is a chance to obtain more refined geometric information.  To begin, one
can directly construct a de Rham representative $C_p(P)$ of $c_p(P)$ from a connection on $P$; this is often also called Chern--Weil theory.  The advantage of the geometric approach is that one can in theory, and sometimes in practice, explicitly compute this de Rham representative from knowledge of the curvature of the connection.

If $C_p(P)$ vanishes pointwise, a very strong condition, then there is a secondary or
Chern--Simons class $TC_p(P)\in H^*(B,\C/\Z)$.  As opposed to the topologically defined 
Chern classes, the Chern--Simons classes depend on the choice of connection, and so are inherently geometric objects.  When defined, the Chern--Simons classes are obstructions to a trivial bundle admitting a trivialization by flat sections of a connection.  Thus the secondary classes are more subtle and correspondingly harder to work with than 
with the primary/Chern classes.  They notably appear as the generators of the integer cohomology of 
the classical groups.  
\medskip

Infinite-dimensional manifolds such as loop spaces $LM$ of manifolds and  mapping spaces 
$\mapsnm$ between manifolds occur frequently in mathematical physics.  Here the structure group of the tangent or frame bundle is an
infinite-dimensional Lie group, the gauge group of a finite rank bundle.  
More generally, a finite rank bundle $E\to M$ over the total space of a fibration 
$M\to B$ of manifolds, the setup of the Atiyah--Singer families index theorem, naturally leads to an infinite rank bundle $\calE\to B$  with a more complicated structure group.
In light of 
physicists' intriguing formal manipulations with path integrals, in particular their quick non-rigorous proofs of the Atiyah--Singer index theorem using loop spaces \cite{atiyah}, it is natural to look for a good theory of characteristic classes of infinite rank 
vector bundles.  

There are several immediate pitfalls.  The fiber of such a bundle, an  infinite-dimensional vector space, comes with many inequivalent norm topologies, in contrast to finite-dimensional vector spaces.  As a result, the topology of the fiber has to be specified carefully.  
If the topology is compatible with  a Hilbert space structure on the model fiber
$\calH$, it is tempting
to take as structure group $GL(\calH)$, the group of bounded invertible operators with bounded inverse.  
However, unlike in finite dimensions, $GL(\calH)$ is contractible, so every $GL(\calH)$ bundle is trivial.  This kills the theory of Chern classes in this generality.  

Of course, $GL(\calH)$ contains many interesting subgroups with nontrivial topology. In
particular, if a subgroup $G$ consists of determinant class operators, then one can try to form the 
characteristic classes associated to the invariant polynomial $\Tr( A^k)$ for $A$ in the Lie algebra of trace class operators;
after all, for finite rank complex bundles,
these polynomials 
form a generating set for the algebra of $U(n)$-invariant polynomials.
However, in infinite dimensions, these classes are generally noncomputable,
in the sense that operator traces are rarely given by e.g., integrals of pointwise calculable expressions.  In particular, it will usually be 
impossible to tell if these Chern classes vanish or not.  Notice that we are not even 
considering the more difficult topological approach of working with $BG$.

In summary, in infinite dimensions we do not expect a version of Chern--Weil theory that applies to all bundles.  Instead, we should look for naturally occurring structure groups with nontrivial topology, and we should look for computable Ad-invariant functions.  

As with the example $\Tr(A^k)$, traces on the Lie algebra of a group give rise to invariant
functions.  For gauge groups, a wide class of traces is known and these ``tend to be"
computable.  However, the determination of all invariant functions is open. This gives us a theory of characteristic classes on mapping spaces, the subject of Section 3, and allows us to 
determine some nontrivial cohomology of mapping spaces.  

This theory is not as geometric as desired, in the sense that  
natural connections on mapping spaces are not compatible with a gauge group, but instead 
are $\pdos$-connections for a larger group $\pdos$ of pseudodifferential operators.  
This larger group has fewer traces, which in fact have been classified.  (An excellent reference for pseudodifferential operators and traces is \cite{scott}.) Again, it is not known
if there are Ad-invariant functions not arising from traces.  

It turns out that the Pontrjagin classes vanish for $\mapsnm$, so we are forced to consider secondary classes.
  In Section 4 we discuss Chern--Simons classes for loop spaces. We use these classes to show that $\pi_1({\diff}(S^2\times S^3))$ is infinite, where ${\diff}(S^2\times S^3)$ is the diffeomrophism group of this 5-manifold.  This result is new but not unexpected, and is given more as an illustration of
potential applications of these techniques.  

In Section 5 we discuss characteristic classes associated to $\diff(Z)$, the group of diffeomorphisms of a closed manifold $Z$.  As pointed out by Singer, there is no known theory of characteristic classes for ${\diff}(Z)$-bundles. Specifically, there are no known
nontrivial Ad-invariant functions on Lie$({\diff}(Z)).$  Instead, we outline a method to detect elements of $H^*({\diff}(G),\C)$ for classical Lie groups $G$.  This is cheating somewhat,
as in finite-dimensional bundle theory we want the cohomology of classifying spaces
like  $BU(n)$, not of $U(n)$ itself.  Of course, $H^*(BU(n),\C)$  
is related to $H^*(U(n),\C)$ by transgression arguments dating back to Borel.   It is
completely unclear if these arguments can be formulated in infinite dimensions, so
the results of this section are baby steps towards understanding characteristic classes for
diffeomorphism groups.

In Section 6, we discuss the setup of the families index theorem of Atiyah--Singer.  As recognized by Atiyah and Singer and used by Bismut, this theorem can be restated in terms of an
infinite rank superbundle $\calE$.  We discuss constructing a theory of characteristic classes on
these bundles.  Here the structure group $\G$ contains both a gauge group and the group
${\diff}(Z)$.  Having a very large group makes it easier to find Ad-invariant functions in principle, but again we know of no nontrivial invariant functions.  
Nevertheless, we can define characteristic classes of $\calE$ for certain connections
due to Bismut.
We discuss an attempt to construct a proof of  the families index theorem using characteristic classes on $\calE.$  While there are serious gaps in the
argument, it is very intriguing that a semidirect product
$\G\ltimes \pdos$ naturally appears as a structure group.
Thus the work in this last section is in some sense is a culmination of the techniques in the 
previous sections. 
\medskip

The determination of the algebra of invariants for Lie groups is a classical topic with a very 19${}^{\rm th}$
century feel.   In highlighting the obvious, namely the central role of these Lie-theoretic questions in Chern--Weil theory, I'm reminded of  Moliere's M. Jourdain, who discovers that he has been speaking prose all his life without knowing it.  In any case, I hope this article spurs interest in extending this classical theory to  infinite-dimensional settings of current interest in geometry and physics.
\medskip

It is a pleasure to thank Andr\'es Larrain-Hubach, Yoshiaki Maeda, Sylvie Paycha, Simon Scott and Fabi\'an 
Torres-Ardila for many helpful conversations on this subject.

\section{{\bf General comments on Chern--Weil theory}}

Let $G$ be a finite-dimensional Lie group, and let $\calP_G$
be algebra of ${\rm Ad}_G$-invariant polynomials from $\calg = {\rm Lie}(G)$ to $\C.$
In its more abstract form, 
Chern--Weil theory gives a map 
$${\cw}:\calP_G\to H^*(BG,\C).$$  
Since $G$-bundles $E\to B$ are 
classified by elements $f\in [B,BG]$, the set of homotopy classes of maps from $B$ to $BG$, a
polynomial $p\in \calP_G$ gives rise to a characteristic class $c_p(E) = f^*{\cw}(p)\in H^*(B,\C).$

For compact connected groups, the suitably normalized map $\cw$ is a ring isomorphism to $H^*(BG,\Z)$ \cite{dieu} (see Ch. 3, Section 4 for references to Borel's original work), 
\cite{dupont}, with the corresponding characteristic classes called Chern classes for $G= U(n)$ and
Pontrjagin classes for $G = SO(n).$    Since the adjoint action is given by conjugation for 
classical groups, for any $k\in \Z^+$ the polynomials $A\mapsto \Tr(A^k)$ are in $\calP_G$.
For $U(n)$, the corresponding characteristic classes are the k${}^{\rm th}$ components of
the Chern character (up to normalization).  It is a classical result of invariant theory that 
these polynomials generate $\calP_{U(n)}$.  
We note that the k${}^{\rm th}$ Chern class
is given by the trace of the transformation induced by $A$ on $\Lambda^k(\C^n)$; since this
transformation is usually also denoted by $A^k$, it is easy to confuse the two uses of 
$\Tr(A^k).$  

\begin{remark} \label{remark1}
  If a group $G$ is linear, i.e., there is an embedding $i:G\to GL(N,\C)$ for some
$N$ (or equivalently, $G$ admits a finite-dimensional faithful representation), then an 
$\adg$-invariant function on $\calg$ corresponds to a $i(G)$-conjugation invariant functional on 
$di(\calg).$  The functionals $A\mapsto \Tr(A^k)$ certainly work, but there may be other 
invariant functions if $i(G)$ has ``small  enough" image in $GL(N,\C)$.  For example, the
Pfaffian ${\rm Pf}(A)$ is an invariant polynomial for $A\in SO(n)$ which is not 
in the algebra generated by  the 
$\Tr(A^k)$; while $\det(A)$ is in this algebra, the Pfaffian  satisfies $({\rm Pf}(A))^2 = \det(A).$
\end{remark}

This abstract approach to characteristic classes is not very useful in practice, both because $BG$ 
tends to be far from a manifold, and because classifying maps are hard to find.  
However, the classifying space approach certainly is powerful.  For example, if a bundle is trivial, then it is classified by a constant map, and so it immediately follows that its characteristic classes vanish.  
(For $U(n)$, the converse almost holds: if a hermitian vector bundle $E$ has  vanishing Chern classes, then some multiple $kE$ is trivial.  The proof is nontrivial.)

There are 
alternative approaches to constructing e.g., Chern classes, one topological and one geometric.  The topological approach constructs the highest Chern class $c_n(E)$ of a 
rank-$n$ complex vector bundle and then iteratively constructs the lower Chern classes by passing to a flag bundle
\cite{huse, milnor}.  Since our bundles will have infinite rank, it's hard to get started on this approach. 

At this point,  we once and for all pass from principal $G$-bundles to vector bundles with $G$ as structure group, although the former case is somewhat more general.  Thus we are assuming that $G$ is linear, and a $G$-bundle denotes a vector bundle with structure group
$G$.

Since the topological approach seems unpromising, we follow the geometric method.  If $B$ is a paracompact manifold, then it admits a partition of unity, so $G$-bundles over $B$
admit $G$-connections.  (Finite-dimensional manifolds and even Banach manifolds are
paracompact).  $c_p(E)$ is then the de Rham cohomology class $[p(\Omega)]$ of $p(\Omega)$,
where $\Omega$ is the $\calg$-valued curvature two-form of the connection,
and $p(\Omega)$ involves wedging of forms and the Lie bracket in $\calg$ in a natural way.
 The Ad-invariance of $p$ is used
crucially to show that $p(\Omega)$ is closed and that its cohomology class is independent
of the connection. This material is standard, and can be found in e.g., \cite{BGV, chern, rosenberg}.  In summary,

\begin{theorem}[Chern--Weil Theorem] \label{CWtheorem}
 Let $p$ be an Ad${}_G$-invariant $\C$-valued power series on
$\calg$.  Let $E\to B$ be a bundle over a manifold $B$ with structure group $G$, and let $\nabla$
be a $G$-connection with curvature $\Omega.$  Then

(i) $p(\Omega)$ is a closed even degree form on $B$.  

(ii) The de Rham cohomology class $[p(\Omega)]\in H^{*}(B,\C)$ is independent of the choice of $\nabla.$
\end{theorem}

Often, the power series of interest are in fact polynomials of some degree $k< {\rm dim}(B)/2$, in which case $[p(\Omega)]\in H^{2k}(B, \C)$.
In particular, for $U(n)$, the classes on a complex bundle $E$ associated to $(2\pi i)^{-1} \Tr(A^k)$ are denoted by
$c_k(E)$ and are called the k${}^{\rm th}$ {\it Chern classes}; the normalization ensures that they are in fact integral classes.
 The {\it Pontrjagin classes} of a real finite rank bundle $F$ are by definition the Chern classes of
the complexification 
$F\otimes \C$, corresponding to the embedding of $SO(n)$ into $U(n)$.  
On $U(n)$, the most important example of a power series is the exponential function
$e^A$; the corresponding Chern class is called the {\it Chern character}. In index theory,
other power series like the $\hat A$-genus naturally occur, although all these are truncated to polynomials at
 the dimension of the manifold. On infinite rank bundles over infinite-dimensional manifolds, there is no reason to truncate, so the use of power series is more natural.

The advantage to the geometric approach is that the de Rham representative $p(\Omega)$
is pointwise computable on $B$.  For example,  a trivial complex bundle has vanishing Chern classes.  From the geometric
construction, we can conclude more:  a bundle with a nonvanishing Chern class does not
admit a flat connection.  This is a stronger statement precisely because there are nontrivial bundles $E$ with flat connections; these have a discrete structure group.

For noncompact finite-dimensional Lie groups, the situation is not so clean.  For $G = GL(n,\R)^+$, the connected group of orientation-preserving elements of $GL(n,\R)$, $\cw$ is not surjective. For
  $G$ is homotopy equivalent to $SO(n)$,
so the universal Euler class, the element $e\in H^n(BSO(n), \Z)$ corresponding to the Pfaffian in 
$\calP_{SO(n)}$, is also an element of $H^n(BGL(n,\R)^+,\Z).$  However, since the Pfaffian
is not $GL(n,\R)^+$-invariant, $e$ is not in the image of $\cw$ on $\calP_{GL(n,\R)^+}.$


It is clear that parts of geometric Chern--Weil theory carry over to infinite-dimensional Lie groups, 
especially since tricky questions about the topology of these Lie groups can often be
avoided.  For example, let $G$ be the group of invertible transformations of a fixed Hilbert space of 
the form $I + A$, where $A$ is trace class.  Since all such operators are bounded, it is not
hard to show that $G$ is indeed a Lie group, with Lie algebra given by the set of trace class
operators. Certainly the first Chern class $c_1(E) = [\Tr(\Omega)]
\in H^2(B,\C)$ exists for any connection on a $G$-bundle $E\to B$.  Here $\Tr$ refers to the operator trace. However,  this Chern class is not computable except in special cases.  In
particular, we do not expect to be able to tell if this class is nonzero or not.  By restricting $A$
to lie in higher Shatten classes, we can construct higher but similarly noncomputable Chern classes, as discussed in \cite{Seg}.

Making sense of Lie groups of unbounded operators on a Hilbert space is difficult, particularly since the exponential map may have a sparse image.  For finite
rank bundles, one can take a default position by considering  bundles for classical groups as $GL(n,\C)$-bundles.  Since
$GL(n,\C)$ deformation retracts onto $U(n)$, the topological theory of characteristic classes
is the same for the two groups.  In fact, the geometric theory is the same, as the Ad-invariant polynomials are the same for the two groups. Of course, for other linear groups the situation can be more complicated.

By analogy, in infinite dimensions we might begin with $GL(\calH)$, the group of bounded
invertible operators with bounded inverses on a real or complex Hilbert space $\calH$.
As an open subset of the set of bounded endomorphisms of $\calH$, $GL(\calH)$ is 
a Lie group \cite{Omori}. However, this group has trivial Chern--Weil theory.
For by Kuiper's Theorem,  the unitary group $U(\calH)$ and hence
$GL(\calH)$ is contractible in the norm topology.  Thus $BGL(\calH)$ has the 
homotopy type of a point, so all $GL(\calH)$ bundles are trivial.
($U(\calH)$ should not be confused with the group $U(\infty) = \lim_{\to}U(n)$, which has nontrivial
topology by Bott periodicity.)  

This problem of having a large contractible structure group 
 is a key feature of infinite dimensions.  As an example of its annoying presence, we can use $GL(\calH)$ to ``ruin" Chern--Weil theory
 for finite rank bundles.  For example, if we embed $GL(n,\C)$ into $GL(N,\C)$ as an upper 
 left block for $N >n$, then we can extend a rank-$n$ complex bundle $E$
 to a rank-$N$ bundle with
 essentially the same transition functions.  Since this amounts to adding a trivial 
 $(N-n)$-rank
 bundle to $E$, the Chern classes are unchanged.  However, if we embed $GL(n,\C)$ into
 $GL(\calH)$, then the extended bundle and any characteristic classes become trivial. (We
 can cook up a similar example in finite dimensions:  Let $E\to S^1$ be the $\Q$-bundle
 which is trivial over $(0,2\pi)$ and with $(0,q)$ glued to $(2\pi, 2q).$  Then $E$ is nontrivial
 as a $\Q$-bundle, but becomes trivial when extended to an $\R$-bundle.)  
 
From these examples, we see that we should consider infinite rank bundles whose structure
group is a subgroup of $GL(\calH)$ with nontrivial topology.  Fortunately, there are several
well-known infinite-dimensional manifolds with a good Chern--Weil theory of characteristic classes.

\section{Mapping Spaces and their characteristic classes}

\subsection{The topological setup}

Let $N^n, M^m$
 be smooth, oriented, compact manifolds.  Fix $s_0 \gg 0$, and let $\mapsnm$
be the functions $f:N \to M$ of Sobolev class $s_0$ (denoted $f\in H^{s_0}$).  Here we fix covers $\{(U_\alpha, \phi_\alpha)\}, \{(V_\beta, \psi_\beta)\}$ of $N, M$, respectively, and we are imposing that
$\psi_\beta f\phi_\alpha^{-1}:\R^n\to \R^m$ is of Sobolev class $s_0$ for all $\alpha, \beta.$
$\mapsnm$ is a smooth Banach manifold \cite{eells}, and the smooth structure is 
independent of the choice of covers.  

We could work with the space of smooth maps from 
$N$ to $M$ as a Fr\'echet manifold, but this is technically more difficult. In particular, the implicit function theorem, which is used repeatedly in the foundations of manifold theory, is
not guaranteed to hold for Fr\'echet manifolds.  This is a little lazy, as the implicit function theorem holds for
tame Fr\'echet manifolds \cite{Hamilton} such as $\mapsnm$, but we choose to work with 
Sobolev spaces just to keep the notation down.  

The easiest examples of mapping spaces are free  loop spaces $LM$ ($N = S^1$) and in particular free
loop groups $LG$.  The most natural bundles are the
tangent bundles $TLM, TLG$.  Just as in finite dimensions, $TLG$ is canonically trivial, so
we do not expect characteristic classes for loop groups.  

To develop a theory of characteristic classes of $TLM$, we should determine its structure 
group and look for Ad-invariant functions.  A tangent vector in $TLM$ at a loop $\gamma$
should be the infinitesimal information in a family of loops $s\mapsto 
\gamma_s(\theta)$, for
$\theta\in S^1$ and $s\in (-\epsilon, \epsilon).$  The infinitesimal information is
$\{ \dot\gamma(\theta)=(d/ds)|_{s=0}\gamma_s(\theta): \theta\in S^1\}$.  This is a vector field along $\gamma$, i.e., a section of
$\gamma^*TM\to S^1$; the pullback bundle has the effect of distinguishing tangent vectors
$\dot\gamma(\theta_0), \dot\gamma(\theta_1)$ where $\gamma(\theta_0) = \gamma(\theta_1).$ Conversely, given a Riemannian metric on $M$,  the exponential maps $\exp_{\gamma(\theta)}:T_{\gamma(\theta)}M
\to M$ combine to take a vector field along $\gamma$ to a loop.  Taking care of the analytic
details, we get $T_\gamma LM = \Gamma(\gamma^*TM)$, where we take $H^{s_0}$ sections
of $\gamma^*TM.$  Since $M$ is oriented, $\gamma^*TM$ is a trivial rank-$m$ real bundle
over $S^1$ denoted $\calR^m = \calR^m_\gamma$.  The trivialization is not canonical, so
$TLM$ need not be trivial, and we have a hope of constructing
characteristic classes.  

We now show that the structure group of $TLM$ is a group of gauge transformations.  
This structure group is determined by the differentials of the transition functions of
$LM$.  Fix a Riemannian metric on $M$. 
Call $s\in T_\gamma LM$ short if  $\exp_{\gamma(\theta)} s(\theta)$ is inside the cut locus 
of $\gamma(\theta)$ for all $\theta.$  Let $U_\gamma$ be the neighborhood of $\gamma$ in 
$LM$ consisting of the all exponentials along $\gamma$ of short loops.  These neighborhoods 
give an open cover of $LM$.
On $U_{\gamma_0}\cap U_{\gamma_1}$, the transition functions are given by fiberwise invertible nonlinear
maps $\Gamma\calR^m\to\Gamma\calR^m$.  Since $\Gamma\calR^m$ is a vector space, the differentials of the transition maps at a $v\in \calR^n$ can be naturally identified
with invertible linear maps on $\calR^m$ which act fiberwise.  Thus the structure group is
 the group of $H^{s_0}$ bundle
automorphisms of $\calR^m$, i.e., the  gauge group $\calG(\calR^n).$  (Strictly speaking, we should take the gauge group of $\calG(T\calR^n)$, but this has the same homotopy type as $\calG(\calR^n).$)

The general case of $\mapsnm$ is similar.  The path components of $\mapsnm$ are in 
bijection with $[N,M]$.  Pick a path component $X_0$ and  $f:N\to M$ in $X_0.$ Then for
all $g\in X_0$, $T_g\mapsnm \simeq \Gamma(f^*TM\to N)$ noncanonically.  $f^*TM$ need
not be trivial, but as above the structure group on $X_0$ is $\calG(f^*TM).$
For convenience, we always complexify real bundles, so the structure group is $\calG(f^*TM\otimes \C).$   From now on, we often omit the $\otimes\ \C$ term.

In summary, $T\mapsnm$ is a gauge bundle, or $\calG$-bundle for short.

Now that the structure group of  $T\mapsnm$ has been
determined, we look for ${\rm Ad}_\calG$-invariant functions on $\frakg$.  Here $\frakg = {\rm Lie}(\calG) = {\rm End}(f^*TM)$ is the vector space of $H^{s_0-1}$ bundle endomorphisms of $f^*TM.$  Since the Lie group and Lie algebra
act fiberwise on $f^*TM$,  the adjoint action of $\calG$
on $\frakg$ is fiberwise conjugation: ${\rm Ad}(A)(b) = AbA^{-1}$.
For  fixed Riemannian metrics on 
$N$ and $M$, $f^*TM$ inherits an inner product, and we can set 
$$c_k:\calG(f^*TM)\to \C,\ \ c_k(A) = \int_N \tr(A^k)\  \dvol_N.$$
Note that the trace depends on the metric on $M$.  This is clearly Ad-invariant, so we can
define
\begin{equation} c_k(X_0) = [c_k(\Omega)]\in H^{2k}(X_0, \C)
\label{cherngauge}
\end{equation}
for $\Omega$ the curvature of any gauge connection on $T\mapsnm.$  We will usually just write
$c_k(\mapsnm)$\ 
$\in H^{2k}(\mapsnm, \C).$ Some examples of these ``gauge classes" will be computed in Section 3.3.

The reader familiar with characteristic classes may be appalled that we are omitting the usual normalizing constants which in finite dimensions guarantee that Chern classes 
have integral periods.  In infinite dimensions, there is no known topological method of producing integral characteristic classes, so there is no natural normalization.

There are certainly many other Ad-invariant functions.  For any smooth function $h:N\to 
\C$, $c_{k,h}(A) = \int_N h\cdot \tr(A^k)\ \dvol_N$ is Ad-invariant.  Letting $h$ approach a 
delta function, we see that for every distribution $h\in \calD(N)$, $h(\tr(A^k))$ is Ad-invariant.
Although it is overkill, this fits in with the finite-dimensional situation, where the structure group $GL(n,\C)$ is the
gauge group of the bundle $\C^n\to *$, where $*$ is a point, and $h\in \calD(*)$ 
must be multiplication by a constant.  
\medskip

\noindent {\it Open question:}  Determine all ${\rm Ad}_\calG$-invariant analytic functions on $
\frakg.$
\medskip

Since $\calG$ is dense in $\frakg$, solving this question includes finding 
 all the traces on $\frakg$, i.e., linear functions $t:\frakg\to \C$ with
$t(ab) = t(ba)$ for all $a,b\in \frakg.$  This is in turn equivalent to computing the Hochschild
cohomology group $HH^0(\frakg, \C)$, which should be feasible. This is interesting even in 
the loop group case, where we are asking for $HH^0(L\frakg, \C)$, where $\frakg$ is now the Lie
algebra of the compact group $G$.  

For an overview of a large class of Ad-invariant functions on $L\calg$ with 
applications to integrable systems, see \cite{sts}.

\subsection{The geometric setup}

Since we are taking a geometric approach to characteristic classes, we should see if the 
natural geometry on $\mapsnm$ is compatible with the structure group $\calG$.  In fact, it is
not, as we now explain.

For simplicity, we will just consider loop spaces $LM^m$.  The parameter $\theta$ always 
denotes the loop parameter, so $LM = \{\gamma(\theta): \theta\in S^1, \gamma(\theta)
\in M\}$ for $\gamma$ of Sobolev class $s_0.$ 

Fix a Riemannian metric on $M$ and fix $s\in [0,s_0]$. Define
an $H^s$ inner product on $T_\gamma LM$ by
$$
\langle X, Y\rangle_{\gamma,s} = \int_{S^1} \langle X_{\gamma(\theta)}, (I+\Delta)^{s}Y_
{\gamma(\theta)}\rangle_{\gamma(\theta)} \ d\theta.
$$
Here $\Delta = D^*D$, with $D = \frac{D}{d\gamma} = \gamma^*\nabla^{M}$ the covariant
derivative along $\gamma$, 
or equivalently the $\gamma$-pullback of the Levi-Civita connection $\nabla^M$ on $M$.
  The role of the
positive elliptic operator $(I+\Delta)^s$ is to count roughly $s-(m/2)$ derivatives of the 
vector fields, by the so-called basic elliptic estimate.  Thus the larger the $s$ and $s_0$, the closer we are to modeling the smooth
loop space.  In particular, the $L^2$ metric (i.e., $s=0$), while independent of a choice of $s$ and hence natural, is too weak for many
situations.  (For example, the absolute version of the Chern--Simons classes discussed in Section 4 are multiples of $s$, and hence vanish at $s=0.$)

This $H^s$ metric gives rise to a Levi-Civita connection $\nabla^{s}$  on $LM$
by the Koszul
formula
\begin{eqnarray}\label{two} 
2\langle\nabla^s_YX,Z\rangle_s &=&\ X\langle Y,Z\rangle_s +Y\langle X,Z\rangle_s
-Z\langle X,Y\rangle_s\\
&&\qquad 
+\langle [X,Y],Z\rangle_s+\langle [Z,X],Y\rangle_s -\langle[Y,Z],X\rangle_s,\nonumber
\end{eqnarray}
but only if the right-hand side is a continuous linear functional of $Z\in T_\gamma LM$.
Note that this continuity is not an issue in finite dimensions.  

\medskip

We will  consider loop groups as an example.  First,  recall that for a finite-dimensional Lie group $G$ with a left-invariant metric, there is a global frame for $TG$ 
consisting of left-invariant vector fields $X_i$.  For such vector fields, 
 the Kozul formula simplifies, since the first three terms
on the right-hand side of (\ref{two}) vanish.  In particular, one can determine $\nabla_{X_j}
X_i$ in terms of the structure constants of $G$.
Since any vector field on $G$ can be written
as $X = f^iX_i$ for $f_i\in C^\infty(G)$,  the Leibniz rule then determines $\nabla_YX$
completely for any $X, Y.$  

For loop groups $LG$, we can take an infinite basis $\{X_i\}$ of $T_{e}LG = L\calg$, by e.g., taking
Fourier modes with respect to a chosen basis of $\calg.$  We can extend these to 
left-invariant vector fields.  A calculation first due to Freed \cite{Freed} gives 
\begin{equation}\label{three}
2\nabla^s_XY=[X,Y]+ (I+\Delta)^{-s}[(I+\Delta)^{s}X,Y]
+(I+\Delta)^{-s}[X,(I+\Delta)^{s}Y],
\end{equation}
for $X, Y$ left-invariant.  The reader is encouraged to rework
 this calculation, which just uses 
that $(I + \Delta)^s$ is selfadjoint for the $L^2$ inner product.  There are technical issues here, such as checking that the 
right-hand side of (\ref{three}) stays in 
$H^{s_0}$, and that   applying the Leibniz rule to infinite sums
$f^iX_i$ also stays in $H^{s_0}.$  Since the $X_i$ are so explicit, these issues can be 
resolved.

In (\ref{three}), $(I+\Delta)^s$ is a differential operator if $s\in \Z^+$ and is a classical
pseudodifferential ($\pdo$) operator otherwise.  In any case, $(I +\Delta)^{-s}$ is always
pseudodifferential.  The critical Sobolev dimension for $LG$ is $1/2$, since loops need to
be in $H^{(1/2)+\epsilon}$ to be continuous, so we will always assume $s> 1/2.$
As an operator on $Y$ for fixed $X$, $\nabla^s_XY$ has order zero:
the first order differentiations in the first and third terms on the right-hand side of (\ref{three})
cancel (as seen by a symbol calculation), and the second term has order $-2s+1$.

These technical calculations are really quite crucial.  On general principles, the connection 
one-form and the curvature two-form take values in the Lie algebra of the structure group. 
So calculating these forms tells us for which structure group $G$ our connection is a $G$-connection. 
For the curvature two-form
$$ \Omega^s(X,Y)=\nabla^s_X\nabla^s_Y-\nabla^s_Y\nabla^s_X-\nabla^s_{[X,Y]},$$
we can always say that
 $\Omega^s\in \Lambda^2(LG, {\rm End}(T_\gamma LG))$ by default, but the vector space
${\rm End}(T_\gamma LG)$ is too big to be useful.  After all,  ${\rm End}(T_\gamma LG)$ could only
be the Lie algebra of ${\rm Aut}(T_\gamma LG)$, all technical issues aside.  Without further restrictions, this group contains both bounded and unbounded operators,
so its topology is unclear.  

Thus without some detailed computations, the setup would be too formal.
However, since $\nabla^s$ is built from zero order $\pdo$s, $\Omega^s$ 
also takes values in zero order $\pdo$s.
That is good news, since 
order zero $\pdo$ are bounded operators on $T_\gamma LG$ with any $H^s$ norm.  
Moreover, the vector space $\pdoo$ of classical $\pdo$s of integer order 
at most zero is the Lie algebra of $\pdos$, the Lie group of invertible classical zeroth order
$\pdo$s \cite{paychasurvey}.  (Inverses of elements in $\pdos$ are automatically bounded.)  
Note that $\pdos \supset \calG(\calR^n)$, since gauge transformations are (zero-th order)
multiplication operators.

Thus by just working out the connection and curvature, we see that it is natural to extend
the structure group from the gauge group $\calG$, which was good enough for the topology of $LG$, to
$\pdos$, which is needed to incorporate the Levi-Civita connection.

Before leaving the loop group case, we note that Freed proved that  $\Omega^s$ actually takes values in $\pdo$s of order at most $-1.$  By some careful calculations, sharp results
have been obtained:

\begin{proposition} \cite{andres}, \cite{MRT}\ If $G$ is abelian, the curvature two-form
$\Omega^s$ for $LG$
takes values in $\pdo$s of order $-\infty$ for all $s \geq 1.$  If $G$ is nonabelian, the curvature two-form takes values in $\pdo$s of order $-\infty$ for $s=1$ and order $-2$ for
$s>1.$  These results are also valid for the based loop groups $\Omega G.$
\end{proposition}

The case $s=1$ is known to be special: for complex groups $G$, 
$\Omega G$ is a K\"ahler manifold for the $s=1$ metric \cite{Seg}. The proposition again 
singles out this case, and 
applies to all
finite-dimensional Lie groups.   One proof that the curvature 
has order $-\infty$ involves showing that the map $\alpha:L\calg \mapsto L\calg[[\xi^{-1}]]$
(which appears in integrable systems as the space of formal nonpositive integer order $\pdo$s on the trivial bundle $E = S^1\times
\calg\to S^1$) given by $\alpha(X) = \sum_{\ell=0}^\infty \frac{(-1)^\ell}{i^\ell}
(\partial_\theta X )\xi^{-\ell}$ is a Lie algebra homomorphism.  It would be interesting to know
how this representation of $L\calg$ on the $H^s$ sections of $E$ fits into the general 
theory of loop group representations.

\medskip

We now consider (\ref{two}) for general loop spaces $LM$. Now all six terms on the 
right-hand
side contribute unless $M$ is parallelizable.  For $s\in \Z^+$, after some simplifications, we end up with terms that
take $\theta$ derivatives of $Z$ (e.g., $(I+\Delta)^s Z$).  This apparently obstructs the 
right-hand side from being a linear functional in $Z$ in the $H^s$ norm, but we can integrate by
parts $2s$ times over $S^1$   to remove this problem.  As a result, the Levi-Civita connection exists in this case.  In contrast, for $s\not\in \Z^+$, trying to integrate by parts
involves the  infinite symbol asymptotics 
$\sigma(I + \Delta)^s \sim \sum_{k=0}^\infty \sigma_{-k}(I+\Delta)^s$ 
of the $\pdo$ $(I+\Delta)^s$, and so
does not terminate.
  To make a long story short, the Levi-Civita connection does not exist.

\begin{theorem} \cite{MRT} Let $M$ be a Riemannian manifold.  Then the Levi-Civita connection 
for the $H^s$ metric on $LM$ exists for $s\in \Z^+ \cup \{0\}.$  If $M$ is not parallelizable,
then the Levi-Civita connection does not exist for $s\not\in \Z^+\cup \{0\}.$
\end{theorem}  

This demonstrates the perils of geometry in infinite dimensions. A similar result should hold for $\mapsnm$ but  has not been worked out.

\begin{remark}
Now that  the structure group $\pdos = \pdos(f^*TM\to N)$ has naturally appeared, we can formulate the notion of 
principal $\pdos$-bundles and associated vector bundles for general 
$\pdos(E\to N)$.  To see that the theory of $\pdos$-bundles is topologically distinct from the
theory of $\calG$-bundles, 
we should show that $\calG$ is not a deformation retract of $\pdos$.  
This is probably true, as for $A\in\pdos$, the top order symbol
$\sigma_0(A)$ lives in $\calG(\pi^*E\to S^*N)$, which
does not retract onto $\calG(E\to N).$  It would be good to work this out.
\end{remark}

To summarize, the geometry of $\mapsnm$ leads us to extend the structure group from a gauge group to a group of $\pdos$.  In contrast, putting a metric on a complex finite rank bundle leads to a reduction of the structure group $GL(n,\C)$ to
$U(n)$.

\subsection{Characteristic classes for $\mapsnm$}

The immediate issue is to find all Ad-invariant functions for $\pdoo = \pdoo(f^*TM\to N).$
 Since this algebra is larger 
than the gauge algebra
$\calg = {\rm End}(f^*TM)$, we expect fewer invariants.  As a start, we know we can build
invariants from traces.  In this setting the traces have been classified, as we explain.

Recall that the Wodzicki residue of a $\pdo$ $A$ acting on sections of a bundle
$E\to N^n$  is defined by
\begin{equation}\label{four} \resw:\pdoo\to\BbC,\ \resw(A) = 
(2\pi)^{-n}\int_{S^*N}
 \tr\ \sigma_{-n}(A)(x, \xi)
\ d\xi\ \dvol(x),
\end{equation}
where $S^*N$ is the unit cosphere bundle.
For dim$(N)>1$, the Wodzicki residue is the unique trace on the full algebra of $\pdo$s up to scaling, although 
the facts that it is a trace and is unique are not obvious \cite{fgl, scott}.  (The issue for
$N= S^1$, the loop space case, is that the unit cosphere bundle is not connected, but this
case has been treated in \cite{ponge}.)  In particular, the integrand 
$ \tr\ \sigma_{-n}(A)(x, \xi)$ is not a trace, so we cannot apply distributions to the integrand to get other traces.
 The Wodzicki residue vanishes on $\pdo$s which do not have a symbol term of order $-n$, so it vanishes on 
$\calg$, on all differential operators, and on all classical operators of noninteger order.  The Wodzicki residue is orthogonal to the operator
trace, in the sense that operators of order less than $-n$ are trace class but have vanishing Wodzicki residue. Just
to reassure ourselves that this trace is nontrivial, for any first order elliptic operator $D$ on sections of $E$,
$\sigma_{-n}(I+D^*D)^{-n/2}(x, \xi) = |\xi|^{-n}$, so $\resw (I+D^*D)^{-n/2}
= {\rm vol}(S^*N)\neq 0.$

\begin{remark}  The Wodzicki residue is the higher dimensional analogue of the residue considered by Adler, van Moerbeke and others in the study of the KdV equation and flows on coadjoint
orbits of loop groups \cite{adler, guest}.  The Ad-invariant functions becomes integrals of motion, and are used
to study the complete integrability of this system. 
\end{remark}

On the subalgebra $\pdoo$, there are more traces.  The leading order symbol trace is 
defined by
\begin{equation}\label{five}
\Tr^{\rm lo}(A) =  
(2\pi)^{-n}\int_{S^*N}
 \tr\ \sigma_{0}(A)(x, \xi)
\ d\xi\ \dvol(x).
\end{equation}
Since $\sigma_0(AB) = \sigma_0(A)\sigma_0(B)$ for $A, B\in \pdoo$, the integrand is a trace, and so any distribution on $S^*N$ applied to the function $\sigma_0(A)\in C^\infty
(S^*N)$ is a trace.

\begin{theorem} \cite{L-P}  For dim$(N) >1,$ all traces on $\pdoo$
are of the form 
$$A\mapsto c\cdot \resw(A) + C(\tr\sigma_0(A))$$
 for some $c\in \C$ and
$C\in \calD(S^*N)$.
\end{theorem}

The proof is an impressive calculation in Hochschild cohomology.  

\begin{remark}   For $\alpha <0$, the set of $\pdo$s of order at most $\alpha$ is a subalgebra
of the full $\pdo$ algebra.  In \cite{L-N}, the traces on these subalgebras are
classified.
\end{remark}

Now that we know what the traces are, we can define two types of characteristic classes
for $\pdos$-bundles for $\pdos = \pdos(F\to Z)$ with $Z$ closed and $F$ a complex bundle. 

\begin{definition}  The k${}^{\rm th}$ Wodzicki--Chern class $\ckw(\calE)$ of the $\pdos$-bundle $\calE\to \calM$ admitting a $\pdos$-connection $\nabla$ is 
the de Rham cohomology class $[\resw(\Omega^k)]\in H^{2k}(\calM, \C)$, where $\Omega$ is the curvature of $\nabla.$  The k${}^{\rm th}$  leading order symbol class $\cklo(\calE)$ is the de Rham class
$[\Tr^{\rm lo}(\Omega^k)].$
\end{definition}

Note that if $\nabla$ is in fact a gauge connection, then $\cklo(\calE)$ is a multiple of the
gauge classes in (\ref{cherngauge}), since for gauge connections the symbol is independent
of the cotangent variable $\xi.$

For $\mapsnm$, we easily get $\ckw(\mapsnm) 
\stackrel{\rm def}{=} \ckw(T\mapsnm\otimes \C) = 0$.  For $T\mapsnm$ is a gauge bundle admitting a  gauge connection.  The curvature form of this connection takes values
in ${\rm End}(f^*TM)$ and so has vanishing Wodzicki residue.  

In contrast, it is shown in
\cite{lrst} that $\cklo(\mapsnm) = {\rm vol}(S^*N)\cdot \ev_n^*c_k(TM)$, where 
the evaluation map $\ev_n:
\mapsnm\to M$ is given by $\ev_n(f) = f(n)$ for a fixed $n\in N.$  With some work, this can be extended to:

\begin{proposition} \label{propprobe}  Let $\mapsfnm$ denote the connected component of 
an element $f\in\mapsnm$ for $M$ connected.
 Let $F \to M$ be a  rank-$\ell$ complex bundle with $c_k(F)\neq 0.$   Then
$$0\neq \cklo(\pi_*\ev^*F)\in H^{2k}(\mapsfnm, \BbC).$$
\end{proposition}


Thus we can use the leading order symbol classes to show that $\mapsnm$ has roughly 
as much cohomology as $M$ does.  Of course, $\mapsnm$ should have much more cohomology.
Setting $M = BU(\ell)$, we get

\begin{theorem}\label{last theorem} \cite{lrst}
Let $E\to N$ be a  rank-$\ell$ hermitian bundle.  There are surjective ring homomorphisms
from $H^*(B\calG(E),\C)$ and $H^*(B\pdos(E),\C)$ to
  the polynomial algebra $H^*(BU(\ell), \BbC) = \BbC[c_1(EU(\ell)),\ldots, c_\ell(EU(\ell))].$
\end{theorem}

This is to our knowledge the first (incomplete) calculation of the cohomology of $B\pdos$, which is
needed for a full understanding of the theory of characteristic classes.  Even for the gauge group, these results seem to be new, although more precise results are known for specific 4-manifolds $N^4$ of interest in Donaldson theory. In contrast, the homotopy groups of (a certain stabilization of) $\pdos$ and hence of $B\pdos$ have been completely computed in \cite{rochon}.  
\medskip

In summary, the study of traces on $\pdoo$ yields two types of characteristic classes, the Wodzicki--Chern classes and the leading order symbol classes, but only the latter are nontrivial.  Note that  we have not addressed
the question of finding Ad-invariant functions not associated to traces.

\medskip

\noindent {\it Open Question:}  Determine all ${\rm Ad}_{\pdos}$-invariant functions on 
$\pdoo.$

\section{Secondary classes on $\pdos$-bundles}

In this section we discuss secondary or Chern--Simons classes in infinite dimensions. This material is taken from \cite{lrst2, lrst,MRT}.

It is useful to think of the Wodzicki--Chern classes as purely infinite-dimensional constructions: if $\calE\to \calM$ is a $\pdos(E^\ell\to *)$-bundle with $*$ just a point,
then $\calE$ is a finite-rank bundle and the only ``$\pdo$s" are elements of $GL(\ell,
\C)$, so there is no Wodzicki residue.  In contrast, the leading order symbol Chern classes 
reduce to the usual Chern classes in this case.

We have already seen applications of the leading order symbol Chern classes.  The Wodzicki--Chern 
classes are poised to detect the difference between $\pdos$- and $\calG$-bundles.  For
as with $\mapsnm$, $\ckw(\calE) = 0$ if $\calE$ admits a reduction to a $\calG$-bundle.
Thus if we can find a single $\pdos$-bundle $\calE$ with $\ckw(\calE) \neq 0$ for some $k$, 
then $\pdos$ cannot have a deformation retraction to $\calG.$  

However, this approach has completely failed to date.

\begin{conjecture}
  For any $\pdos$-bundle $\calE$ over
a paracompact base, $\ckw(\calE) = 0$ for all $k$.
\end{conjecture}

This conjecture holds if either (i) the structure group of $\calE$ reduces from $\Psi_0^*$ to
the group ${\rm Ell}^*$ of invertible zeroth order $\pdo$s with leading symbol the identity \cite{lrst}, or
(ii) $\calE$ admits a bundle map via fiberwise Fredholm zeroth order 
$\pdo$s to a trivial bundle \cite{lrst2}.

The proof of (i) uses the fact that ${\rm Ell}^*$ has the homotopy type of invertible operators of the form identity plus smoothing operator, and these operators have vanishing Wodzicki residue.  

For (ii), we first note that this condition always holds for finite-rank bundles.
The proof follows the structure of the heat equation proof of the families index theorem (FIT) for superbundles \cite{Bismut2}.  One takes a superconnection $\nabla$ on $\calE$ and modifies it to
$B_t = \nabla + t^{1/2}A$, where $A$ is a zero-th order odd operator on each fiber.  The Wodzicki--Chern
character of $\calE$ has representative $\exp(-B_t^2)$ for any $t \geq 0.$  As $t\to\infty$, 
$B_t$ becomes concentrated on the finite rank index bundle for $A^*A$, and so
the Wodzicki--Chern character vanishes.   This implies that all the Wodzicki--Chern classes
vanish.  

\begin{remark}\label{althoughalmostall} Although almost all details have been omitted, this proof is much quicker than the heat equation proofs of the FIT, and for good reason.  In the FIT proofs, one is trying to compute the Chern character of the index bundle in terms of characteristic classes by 
comparing the $t\to 0$ and $t\to \infty$ limits of heat operators.  The $t\to\infty$ limit is
relatively easy, and is mimicked in the proof outlined above.  However, because one essentially wants to use the operator trace, the construction of the appropriate $B_t$ is 
much more delicate, in order to have a well-defined limit at $t=0$.  Once again, we see
the strong contrast between the operator trace and the Wodzicki residue.
\end{remark}

The vanishing of the Wodzicki--Chern classes is not the end of the story, as it is in fact
the prerequisite to defining secondary classes.  Let
 $\nabla_0, \nabla_1$ be connections on a $\pdos$-bundle with local connection one-forms
 $\omega_0, \omega_1$ and curvature $\Omega_0, \Omega_1$.
Then just as for
finite-rank bundles, as even forms
\begin{equation}\label{4.0}\ckw(\Omega_1) - \ckw(\Omega_0) = 
d\csw_k(\nabla_1,\nabla_0),
\end{equation}
 with the odd
 form $\csw_k(\nabla_1, \nabla_0)$ given by
\begin{equation}\label{5.1}
\csw_k(\nabla_1,\nabla_0) = \int_0^1 \resw[(\omega_1-\omega_0)\wedge \overbrace{\Omega_t\wedge ...\wedge \Omega_t}^{k-1}]
\ dt,
\end{equation}
where 
$$\omega_t = t\omega_0+(1-t)\omega_1,\ \ \Omega_t = d\omega_t+\omega_t\wedge\omega_t.$$ 
Here we are just lifting the finite-dimensional formula from \cite[Appendix]{chern}, replacing
the matrix trace with the Wodzicki residue.
(\ref{4.0}) is precisely the explicit formula showing that the Wodzicki--Chern class is independent of
connection, and so is the $\pdos$ version of the proof of Theorem \ref{CWtheorem} (ii).  

(\ref{4.0}) shows that $\csw_k$ determines a $2k-1$ cohomology class if the 
Wodzicki--Chern forms for $\nabla_0, \nabla_1$ vanish pointwise.  This holds in all
cases we have been
 able to compute.

\medskip

\noindent {\it Open question:}  Do Wodzicki--Chern forms always vanish pointwise?
\medskip

If this is the case, then the theory of secondary classes for $\pdos$-bundles based on the Wodzicki residue 
produces cohomology classes in odd degrees.

For finite-rank bundles $E\to N$, the corresponding classes are called Chern--Simons classes
$CS_k(\nabla_0, \nabla_1)\in H^{2k-1}(N,\C)$, when they exist.  
In contrast to the Chern
classes, which are defined via the geometric Chern--Weil theory but have topological content, 
these ``relative"  Chern--Simons classes really do depend on the choice of two
connections, and so are 
geometric objects.  
There is an ``absolute" version of Chern--Simons classes $CS_k(\nabla)\in H^{2k-1}(N, 
\C/\Z)$ that only uses one connection but takes values in a weaker coefficient ring
\cite{C-S}.  

\begin{definition}  The k${}^{\rm th}$ Wodzicki--Chern--Simons (WCS) class associated to 
connections $\nabla_0, \nabla_1$ on a $\pdos$-bundle $\calE\to \calM$
is the 
cohomology class of $\csw_k(\nabla_1,\nabla_0)$ in $H^{2k-1}(\calM, \C)$, provided this form is closed.
\end{definition}

In finite dimensions there are two ways to assure that characteristic forms for 
$E\to N$ vanish:  either
use a flat connection, or pick a form whose degree is larger than the dimension of $N$
or the rank of $E$.  For
example, if $E$ is trivial, it admits a flat connection $\nabla$, so we can define
$CS_k(\nabla, g^{-1}\nabla g)$ for any gauge transformation $g$.  The dimension restriction
is only useful to define $CS_{r}$ for dim$(N)= 2r-1$; this was used very effectively by
Chern--Simons \cite{C-S} and Witten \cite{witten} to produce invariants of 3-manifolds. 

For $\mapsnm$ with $M$ parallelizable, a fixed trivialization of $TM$ leads to a global trivialization of $T\mapsnm$. A gauge transformation $g$ of $TM$ induces a gauge transformation of $T\mapsnm$, so one has an element
$\csw_k(\nabla, g^{-1}\nabla g)$
$\in H^{2k-1}(\mapsnm, \C).$  

\medskip

\noindent {\it Open question:}  Is $\csw_k(\nabla, g^{-1}\nabla g)$ ever nonzero?
\medskip

Although the dimension restriction on $N$ looks incapable of generalization to\\
 $\mapsnm$, 
an examination of the representative $\resw((\Omega^s)^k)$ of $\ckw$ for the $H^s$ metric
($s\in\N$) shows that this form
vanishes pointwise for $2k > {\rm dim}(M)$, as in (\ref{5.4}) below.  Thus we always get a secondary class
$\csw_k(\nabla_1, \nabla_0)\in H^{2k-1}({\rm Maps}(N, M^{2k-1}), \C)$   associated to the 
$s=1, 0$ metrics on $\mapsnm$ determined by fixed metrics on $N, M$.  

For simplicity, we go back to the loop space case $LM$. Because the 
$L^2$ (or $s=0$) connection is so easy to treat -- its connection one-form is just the one-form
for the Levi-Civita connection on $M$ -- the local formula for $
\csw_k = \csw_k(\nabla_1, \nabla_0)$
really is explicitly computable \cite[Prop. 5.4]{MRT}:  as a $(2k-1)$- form on $T_\gamma LM$, 
\begin{eqnarray}\label{5.4}
\lefteqn{CS^{\rm w}_{k}(X_1,...,X_{2k-1}) }\\
&=&
\frac{2}{(2k-1)!} \sum_{\sigma} {\rm sgn}(\sigma) \int_{S^1}\tr[
(-\Omega^M(X_{\sigma(1)},\dg)
-2\Omega^M(\cdot,\dg)X_{\sigma(1)}) \nonumber\\
&&\qquad 
\cdot\  (\Omega^M)^{k-1}(X_{\sigma(2)},\ldots,X_{\sigma(2k-1)} )],\nonumber
\end{eqnarray}
where $\sigma$ is a permutation of $\{1,\ldots,2k-1\}$ and $\Omega^M$ is the curvature of the
Levi-Civita connection on $M$.  The tangent vectors $X_i$ are vector fields in $M$ along
$\gamma$, so we see that this form will vanish if $2k-1> {\rm dim}(M).$

We would like to use $\cswk$ to detect odd cohomology in $H^*(LM,\C).$  We need a test
cycle in degree $2k-1$.  The natural candidate is $M$ itself, thought of as the set of constant loops.
However, for these loops $\dg = 0$, so (\ref{5.4}) vanishes.  Instead, assume that
$M$ admits an $S^1$-action $a:S^1\times M\to M.$  This induces a map $
\tilde a:M\to LM$ by $\tilde a(m)(\theta) = a(\theta, m).$  Dropping the tilde, the action now
produces a test cycle $a_*[M]\in H_{2k-1}(LM, \Z)$, where $[M]$ is the fundamental class of 
$M$, which is assumed orientable.  (For the trivial action, $a_*[M]$ is $M$ as constant loops.)  

If $\int_{a_*[M]} \cswk = \int_M a^*\cswk$ is nonzero, then $\cswk\neq 0$ in the cohomology of $LM$.  The computation of the integral is frustrating:  it always vanishes in
the easiest case ${\rm dim}(M) = 3$; in higher dimensions, 
most explicit Riemannian metrics come with large continuous symmetry groups, and in all
examples we get
$ \int_{a_*[M]} \cswk = 0$,  although we cannot prove a general vanishing theorem.  Fortunately, there is a  family $g_t, t\in (0,1),$ of Sasaki--Einstein metrics on 
$S^2\times S^3$ due to \cite{gdsw} which is explicit enough and
has enough symmetry to make the integral
calculation feasible but which is not ``too symmetric."  In particular, this construction gives a metric fibration $S^1\to S^2\times S^3\to S^2\times S^2$ generalizing the Hopf fibration, so we get
a circle action by rotating the fiber.
The calculations can be done in closed form by
Mathematica${}^\copyright$.  We get $\int_{a_*[M]} \cswk\neq 0$, and so with a little work we conclude that  $H^5(L(S^2\times S^3), \C)$ is infinite.  

A circle action $a$ on $M$ also induces $\tilde a: S^1\to {{\diff}}(M)$ by $\tilde a(\theta)(m) = a(\theta, m)$, 
so we get an element of $\pi_1({{\diff}}(M)).$
It is easy to check that $\int_{a_*[M]} \cswk \neq 0$ implies $\pi_1({{\diff}}(M))$ is infinite.  In
particular, $\pi_1({{\diff}}(S^2\times S^3))$ is infinite.

We tend to trust
the computer calculations, as in the $t\to 0$ limit the metric $g_t$ becomes a metric on 
$S^5$ and the integral explicitly vanishes.  This matches with the known result that $\pi_1({{\diff}}(S^5))$ is finite.  On the other hand, up to factors of $\pi$, the integrals calculated are
always rational; this needs further explanation.

\begin{remark}  (i) While there is nothing in theory to stop us from computing WCS classes for
$\mapsnm$, in practice the number of computations necessary to compute the Wodzicki
residue of an operator on $N^n$ increases exponentially in $n$.  So while computations
are feasible for loop spaces and ${\rm Maps}(\Sigma^2, M)$, the setting for string theory, 
one needs a very good reason to do computations on higher-dimensional source manifolds.

(ii) These results are too specific.  If one can show that the Wodzicki--Chern forms always vanish pointwise, then a much more robust theory of WCS classes would be available.
\end{remark}

\section{Characteristic classes for diffeomorphism groups}

The search for characteristic classes associated to the diffeomorphism group of a
closed manifold $X$ can be interpreted in two ways:  (i)  ${{\diff}}(X)$ is an open subset of 
${\rm Maps}(X,X)$, and so as in Section 3
characteristic classes can be used to detect elements of 
$H^*({{\diff}}(X), \C)$;  (ii) certain infinite rank bundles have ${{\diff}}(X)$ as 
part of their structure group.  In this section we consider the first question, and in Section 6 we treat (ii).

First, the proof that cohomology classes for ${\rm Maps}(X,X)$ 
are detected by characteristic classes of $X$ (Prop.~\ref{propprobe}) 
 does not carry over to ${{\diff}}(X)$.  
Indeed, the proof should break down, since as in finite dimensions the Lie group ${{\diff}}(X)$ 
admits a flat connection and so has vanishing leading order symbol classes.  (As usual,
there are technicalities about the Lie group structure on ${{\diff}}(X)$, which are most easily 
treated by considering $H^s$ diffeomorphisms.)


Thus we expect to find only secondary classes.  
We now outline a method that may produces odd degree classes in $H^*({\diff}(X), \C)$.

The cohomology ring  $H^*(U(n),\Z)$ is generated by suitably normalized
Chern--Simons classes, as the standard forms 
$\tr((g^{-1}dg)^{2k-1})$ built from the Maurer--Cartan form $g^{-1}dg$ are the Chern--Simons forms for $c_k$
associated to the flat connection $\nabla$
on the trivial bundle $U(n) \times {\mathfrak gl}(n, \C)$ and to the gauge equivalent connection $g^{-1}\nabla g$ \cite{mis-r-t}.  Here we think of $g$ as the gauge transformation $M\mapsto g\cdot M$ for $M\in {\mathfrak gl}(n,\C).$
There are similar results for other classical linear groups; see e.g., \cite[Ch. 4.11]{Seg} and the Bourbaki references therein,  \cite{c-e}, particularly the references to the original work of Hopf, and \cite[Ch. 3D]{hatcher} for a modern treatment for $SO(n)$.  Moreover, 
a certain average of these forms along loops in $G$ give generators for $H^*(LG, \Z)$ 
\cite{Seg}.


In finite dimensions, identifying the Maurer--Cartan form with a gauge transformation requires an embedding $G\to GL(N,\C).$  For $G={\diff}(X)$, 
 we
can embed $i:G\to GL(\Gamma(\calC^N))$, where
$\calC^N = X\times \C^N$ is  the trivial bundle over $X$,
  $GL$ refers to bounded 
operators with bounded inverses, and $\Gamma(\calC^N)$ refers to $H^s$ sections. The
 embedding is given by $i(\phi)(s)(x) = s(\phi^{-1}x).$  Now $\phi$ makes sense as a gauge
  transformation of the trivial bundle ${\diff}(X)\times \Gamma(\calC^N)$
   via $s\mapsto i(\phi)(s).$

We can now define 
the Maurer--Cartan form $\phi^{-1}d\phi$ for $\phi\in {\diff}(X)$. As in Section 4, we will get secondary classes in $H^{\rm odd}({\diff}(X), \C)$ associated to the trivial connection $\nabla$ and 
$\phi^{-1}\nabla\phi$
once we pick an ${\rm Ad}_G$-invariant function on 
${\rm Lie}({\diff}(X)) $.  Since a family of diffeomorphisms $\phi_t$ starting at the identity has
infinitesimal information $\dot\phi_0$, a vector field on $X$, we have
${\rm Lie}({\diff}(X))= \Gamma(TX).$  
The adjoint action of a diffeomorphism $\phi$ on a vector field $V$ is easily seen to be $V\mapsto \phi_*V$.  
\medskip

\noindent {\it Open Question:}  Find a nontrivial ${\rm Ad}_{{\diff}(X)}$-invariant function  on the set of ($H^s$) vector fields on a 
closed manifold $X$.
\medskip

\begin{remark} (i) Because the adjoint action is not by conjugation, finding a trace on $\Gamma(TX)$ does not produce secondary classes.  It seems to be an open question whether there exist any nontrivial traces on $\Gamma(TX)$ for general $X$.
%
 If $X=G$ is itself a compact linear Lie group, there are many traces on $T_I{\diff}(G) = 
\Gamma(TG)$.  Namely, a vector field $V$ on $G$ is just a $ \calg$-valued function on $G$, 
so for any distribution $f\in \calD(G)$, $f(\tr(V))$ is a trace. This case needs further work.

(ii)  There is another context in which this Open Question comes up.  
${\diff}(X)$ is the structure group for fibrations $X\to M \to B$ of manifolds with fibers $X$, so characteristic classes associated to ${\diff}(X)$ would be obstructions to the triviality of a fibration, just as ordinary characteristic classes as obstructions to the triviality of principal
$G$-bundles. 
\end{remark}

Thus this general approach to secondary classes for ${\diff}(X)$ is unavailable at present.
However,
for $X=G$ a compact linear Lie group, we can detect odd degree classes in $H^*({\diff}(X),\C)$ using finite rank bundles.  Let $\alpha:G\to{\diff}({G})$ be the embedding $g\mapsto L_g$, for $L_g$ left translation by $g$. 
The trivial rank bundle $\calC_1^N = {\diff}({G})\times {\mathfrak gl}(N,\C)$ admits the gauge transformation
$\phi\in{\diff}({G})\mapsto (M\mapsto \phi(e)M)$.  This gauge transformation, 
also denoted by $\phi$,
 restricts to the gauge transformation $g$ on $\alpha(G)\subset {\diff}(G)$ 
 since $L_g(e) = g.$  The bundle $\calC_1^N$ 
 has the trivial connection $\nabla_1$. The gauge transformed connection 
 $\nabla_1^\phi = \phi^{-1}\nabla_1\phi$ has the global connection one-form $\phi^{-1}d\phi$,
 which restricts to $g^{-1}dg$ on $\alpha(G).$   
On the finite-rank bundle $\alpha(G)$, we
use the ordinary matrix trace to define Chern--Simons forms
$$CS^{2k-1}(\nabla_1,\nabla_1^\phi) = \tr ( (\phi^{-1}d\phi)^{2k-1}).$$
We can also define Chern--Simons classes for $\nabla^g$ on  $\calC^N$ by the same formula.
A straightforward calculation gives
$$ \int_{\alpha(z_{2k-1})} CS^{2k-1}(\nabla_1, \nabla_1^\phi) =
  \int_{z_{2k-1}} CS^{2k-1}(\nabla,\nabla^g)$$
for any $(2k-1)$-chain $z_{2k-1}$ on $G$.  
This implies

 \begin{theorem} \cite{mis-r-t} For any compact linear group $G$, the map $\alpha:G\to{\diff}
     ({G})$, $g\mapsto L_g$, induces a surjection
 $\alpha^*: H^*({\diff}({G}),\R)\to H^*(G,\R).$  
 \end{theorem}
 
 Dualizing this result, we see that the real homology of $G$ injects into the homology of ${\diff}(G)$.  As with mapping spaces, we expect nontrivial homology in infinitely many degrees, but these techniques only give information up to dim$(G)$.

\section{Characteristic classes and the Families Index Theorem}

As explained below, the families index theorem (FIT) is a generalization of the
Atiyah--Singer index theorem.  Infinite rank superbundles $\calE$ naturally appear in the setup of the FIT.  In this section we discuss how a theory of characteristic classes on 
$\calE$ may give insight into the FIT.  In particular, the relevant structure group incorporates aspects of gauge groups, diffeomorphism groups, and the group $\pdos$ of
pseudodifferential operators discussed in previous sections.  This section is based on \cite{lprs}.

The bundle $\calE$ was explicitly mentioned by Atiyah and Singer \cite{AS-IV} as a topological object, but was not used in their proof.  Bismut \cite{Bismut2} used $\calE$ 
as a geometric object, in that he constructed what is now called the Bismut superconnection on
$\calE$.   Bismut did not define characteristic classes for $\calE$, because he did not
 need them:
the fine details of his proof  take place on the sections
of a finite rank bundle, the model fiber of $\calE$ (see Remark \ref{bproof}).  In this section we
try to define characteristic classes directly  on $\calE$.  The hope, not yet realized, is that
not only is $\calE$ a proper setting for the FIT, but that a proof of the FIT can take place 
on $\calE.$

We recall the basic setup, inevitably leaving out a slew of technicalities.  
Let $Z\to M\stackrel{\pi}{\to} B$ be a fibration of closed  connected
manifolds, and let $E, F\to M$ be 
finite rank bundles.  Set $Z_b = \pi^{-1}(b)$, $E_b = E|_{Z_b}, F_b = F|_{Z_b}.$  Assume that
we have a smoothly varying family of elliptic operators $D_b:\Gamma(E_b)\to \Gamma(F_b).$  Although the dimensions of the kernel and cokernel of $D_b$ need not be continuous in $b$, the index ind$(D_b) = {\rm dim\ ker}(D_b) - {\rm dim\ coker}(D_b)$ is 
constant.  It is therefore plausible and indeed true that the virtual bundle IND$(D) = 
[{\rm ker}(D_b)] - [{\rm coker}(D_b)]\in K(B)$ is well defined.  

Although the FIT can be stated entirely within K-theory, it is easier to state it as an 
equality in cohomology.  The Chern character $ch:K(B)\otimes\Q\to H^{\rm ev}(B, \Q)$ is an
isomorphism, and the FIT identifies $ch({\rm IND}(D))$ with an explicit characteristic class 
built from the symbols of the $D_b.$  

Even the case where $B=\{b\}$ is a point is highly nontrivial.   In this case, 
$$ch({\rm IND}(D))= {\rm ind}(D_b)\in H^0(\{b\}, \Q) = \Q$$
 (of course the index is an integer). Identifying the corresponding characteristic class gives the ``ordinary" Atiyah--Singer index theorem, which generalizes Riemann--Roch type theorems 
 for smooth varieties and the Chern--Gauss--Bonnet theorem.  Thus in the 
 appearance of a base parameter space, the FIT is a smooth version of
 Grothendieck--Riemann--Roch theorems.
 
 Rather than discussing the characteristic classes built from the symbols of general $D_b$, we will discuss the particular case of families of
 coupled Dirac operators; in fact, a K-theory argument shows that proving the FIT for coupled Dirac operators implies the full FIT.  (Coupled Dirac operators are discussed in e.g.,
 \cite{l-m}.) 
 So assume that $M$ and every fiber $Z_b$ are orientable and spin in a compatible way, and that $E = S^+\otimes K, F = S^-\otimes K$, where $S^{\pm}$ are the spinor bundles for $M$ and $K$ is yet another bundle over $M$.  Put a metric on $M$; the restriction of the metric to each $Z_b$ defines a Dirac operator on $S^{\pm}_b$.  Put a connection $\nabla^K$ on $K$.  This induces connections $\nabla_b$ on each $K_b$, and gives a family of coupled Dirac operators 
 $\ddd^{\nabla^K} = \ddd^{\nabla_b}: E_b\to F_b.$  
 
 Let $\hat A(M)$ be the $\hat A$-genus of $M$, and let $\int_Z$ denote integration over the fiber 
 (i.e., capping with $Z$ as a class in $H_*(M, \Q)).$
 
 \begin{theorem} (FIT) $ch({\rm IND}(\ddd^\nabla)) = \int_Z \hat A(M)\cup ch(K)$ in $H^{\rm ev}(B, \Q).$
 \end{theorem}
 
 If $\pi_!:K(M)\to K(B)$ is the analytic pushforward map, which by definition sends
 $H$ to $\indd(\ddd^{\nabla^H})$, then the FIT can be restated as the commutativity of the diagram
$$ \begin{CD}  K(M)@>{ch}>> H^{\rm ev}(M,\Q)\\
@V{\pi_!}VV @VV{\int_Z \hat A\ \cup\ (\cdot)} V\\ 
K(B)@>{ch}>> H^{\rm ev}(B, \Q).\\
\end{CD}
 $$
 
As in the last section, the structure group of the fibration is ${\diff}(Z)$, but there is more going on.  The bundle $E\to M$ pushes down to the infinite rank bundle $\pi_*E = \calE
\to B$, where the fiber $\calE_b$ is the smooth sections of $E_b$.  (Thus $\calE$ is the sheaf theoretic  pushdown of the sheaf $\Gamma(E).$)  It is easy to check that 
$\calE_b$ is modeled on $\Gamma(F)$ for some bundle $F
\stackrel{p}{\to} Z$ of rank equal to rank$(E)$.  

  The structure 
group of $\calE$ is 
\begin{equation*}\label{strgp}
\G = \left\{ 
\begin{CD}  F@>{f}>> F\\
@V{p}VV @V{p}VV\\ 
Z@>{\phi}>> Z\\
\end{CD}  : \phi\in {\diff}(Z), f|_{F_{z}}\ {\rm a\ linear\ isomorphism}\right\}.
\end{equation*}
This just says that when fibers of $M$ are glued by a diffeomorphism of $Z$ over $b\in B$,
$E_b$ must be glued by a bundle isomorphism.  
$\G$ is called ${\diff}(Z,F)$ in \cite{AS-IV}.
These transition maps act pointwise on $\Gamma(F)$ (within a fixed Sobolev class), the model space for the fibers of $\calE$,  by 
\begin{equation*}\label{action}
s\mapsto fs\phi^{-1},\ \ {\rm i.e.,}\ \  s\mapsto [z\mapsto f(z)(s(\phi^{-1}(z)))].
\end{equation*}
We note that this equation defines a faithful action of $\G$ on $\Gamma(F)$, which is a Hilbert space once we fix a Sobolev class of sections.
The tangent space to $\G$ at a pair $(\phi, f)$ is given by \cite{Omori}:
$$T_{(\phi, f)} \G = 
\left\{ \begin{CD}  F@>{s}>> f^*TF\\
@V{p}VV@V{f^*p_*}VV\\ 
Z@>{V}>> \phi^*TZ\\
\end{CD}  : s|_{F_{z}} \ {\rm linear}\right\}.
$$
This follows from thinking of $\phi$ as an element of ${\rm Maps}(Z,Z)$ and similarly for $f$, and calculating as in previous sections.
In particular, the Lie algebra  $\calg =T_{(id, id)}\G$ is
\begin{equation*}\label{one}
\calg = 
\left\{ \begin{CD}  F@>{s}>> TF\\
@V{p}VV@V{p_*}VV\\ 
Z@>{V}>> TZ\\
\end{CD} : s|_{F_{z}} \ {\rm linear}\right\}.
\end{equation*}

The difficulty of implementing Chern--Weil theory is the following: 

 \medskip
 \noindent {\it Open question:}  Find a nontrivial ${\rm Ad}_{\G}$-invariant function on $\calg.$
 \medskip
 
 \begin{remark}\label{remark2}  (i) The subgroup of $\G$ where $\phi$ is the identity is precisely the gauge
 group $\calG(F)$ of $F$, so we are looking for generalizations of the invariant functions in Section 3.  
The structure group restricts to $\calG(F)$ if the fibration $Z\to M\to B$ is trivial.  In particular, if the fibration is $N\to \mapsnm\times N \to \mapsnm$ and $E = \ev^*TM\to\mapsnm$, then $\calE$ is precisely $T\mapsnm.$  

(ii) As a subcase of (i), fix 
a compact group $G$,  let $E$ be a trivial $G$-bundle, and let 
$S^1\to M = B\times S^1\to B$ be a trivial circle fibration.  Then $\calG(F)$ is the loop group $LG$.  This
so-called caloron correspondence between $G$-bundles over $M$ and $LG$-bundles over $B$ is discussed thoroughly in \cite{m-v} and goes back to work of Garland and Murray
\cite{g-m}.  In particular, characteristic classes for $LG$-bundles are constructed by
Murray and Vozzo in 
\cite{m-v}; the characteristic classes treated below reduce to the Murray--Vozzo classes
in this case.
\end{remark}

We can  avoid answering the Open  Question and still define characteristic
 classes in a restricted sense.  We first recall the construction of a connection on $\calE$ due to Bismut
 \cite{Bismut2}.   Let 
$HM$
be a complement to the vertical bundle $VM = {\rm ker}\ \pi_*$ in $TM.$  For example, if we have chosen a metric on $M$, we can take $HM =(VM)^\perp.$  
Recall that $E\to M$ has a connection $\nabla^E$; for a given hermitian metric on $E$,
we may assume that $\nabla^E$ is a hermitian connection.  The Bismut
 connection $\nabla = \nabla^B$ on $\calE\to B$ is defined by
\begin{equation}\label{bisconn1}
\nabla_X r(b)(z) = \nabla^E_{X^H} \tilde r(b,z),
\end{equation}
where $X\in T_bB, r\in\Gamma(\calE), z\in \pi^{-1}(b)$, $X^H$ is the horizontal lift of $X$ to
$HM_{(b,z)}$, and $\tilde r\in \Gamma(E)$ is defined by $\tilde r(b,z) = r(b)(z).$  (Here we 
abuse notation a little by writing $(b,z)$, which assumes that a local trivialization of $
Z\to M\to B$ has been given.)

\begin{remark} \label{bproof} We outline Bismut's heat equation proof of the FIT.  Bismut first adjusts
$\nabla^B$  to be unitary with respect to the $L^2$ hermitian metric on $\calE$.  He then 
modifies the new connection in a nontrivial way to form a superconnection $\nabla_t$ on
$\calE$ for $t>0.$  The ``curvature" two-form $\nabla_t^2$ acts fiberwise on $\calE$, just
as for finite rank bundles, and takes values in smoothing (and hence trace class) operators.  As $t\to\infty$, the form-valued operator trace $\Tr(\nabla_t^2)$ converges to a representative of the Chern
character of the index bundle; this step is not too difficult in light of the original heat equation
proof of the index theorem.  
As $t\to 0,$  $\Tr(\nabla_t^2)$ converges nontrivially to the differential form
representative of the right hand side of the FIT.  
It is not hard to show that the 
two  limits differ
by an exact form, so their cohomology classes are the same.  

This is called the local form of the FIT, since the proof generates the specific characteristic forms in the right cohomology class.
\end{remark}   

The (easy) Bismut connection fits into the Atiyah--Singer framework as follows:

\begin{lemma} \cite{lprs} \label{bisconn} The Bismut connection is a $\G$-connection.  In a fixed local trivialization,
the connection one-form assigns to $X\in T_bB$ the pair $(V,s)\in \calg$, where $V = 
\dot\phi_t(0) = X^H$ and $s(v)_{(b,z)} = (d/dt)_{t=0}\Vert_{0,t}(z) v.$
\end{lemma}

Here the parallel translation $\Vert_{0,t}$ for the Bismut connection is thought of as a bundle isomorphism of $F$ via a local trivialization.  The proof directly shows that the 
holonomy of the Bismut connection lies in $\G$.  

This suggests that if we can define the Chern character for connections on $\calE$ 
 for the coupled Dirac operator case (i.e., $\calE$ is the superbundle associated to $(S^+\otimes K)$\\
 $ \oplus\ (S^-\otimes K)$ ), then we could try to 
prove the local FIT by showing

\medskip
(i) For the (easy) Bismut connection on $\calE$, the representative  differential form
$ch(\Omega^B)$
of the Chern character $ch(\calE)$
equals $\int_Z \hat A(\Omega^M) ch(\Omega^K)$, where $\Omega^M$
is the curvature of the Levi-Civita connection for the metric on $M$, and $\Omega^K$ is
the curvature of the connection on $K.$

(ii) There exists a connection on $\calE$ for which $ch(\calE) = ch(\indd(\ddd^{\nabla^K})) \in H^{\rm ev}(B, \Q).$
\medskip

Since the Chern character should be independent of the connection on $\calE$, this would give the FIT.  

As we will now see, step (i) fails, but in a very precise way.

For a fixed hermitian connection $\nabla^E$ on $E$, the
associated Bismut connection has
 curvature two-form $\Omega^B$ taking values in $\calg$, so we can write
$\Omega^B(X,Y) = (V, s)$ for $X, Y\in T_bB$, $V\in \Gamma(TZ)$, and $s\in \Gamma(TF).$  With respect
to a local trivialization, we can consider $s\in\Gamma(TE_b).$
 The connection $\nabla^E$ induces a connection on $E_b$, or equivalently gives a splitting 
$TE_b = VE_b\oplus HE_b.$  The vertical component $(\Omega^B)^v = s^v\in VE_b$ can
naturally be identified with a map $s^v:E_b\to E_b$.  Since $s$ covers $V$, it easily 
follows that $s^v\in {\rm End}(E_b),$ i.e., $s^v$ is a fiberwise endomorphism of $E_b.$  Thus we can create forms 
\begin{equation}\label{cwf}
b\mapsto c_k(\Omega^B_b) = \int_Z\tr ((\Omega_b^B)^{v})^k\in \Lambda^{2k-{\rm dim}(Z)}(B).
\end{equation}

We claim that these forms are closed and have de Rham class independent of the choice of $\nabla$  
on $E$.
To see this, first note that it is 
well known and not 
difficult to compute that the Bismut connection on $\calE$ has curvature
\begin{equation}\label{Bcurv}\Omega^B(\xi_1, \xi_2) = \nabla^E_{T( \xi_1^H, \xi_2^H)}+ R^E(\xi_1^H, \xi_2^H),
\end{equation}
with $T( \xi_1^H, \xi_2^H) = [\xi_1, \xi_2]^H - [\xi_1^H, \xi_2^H]$
and $R^E$ the curvature of $\nabla^E$.  Moreover, the first term
on the right-hand side of (\ref{Bcurv}) is horizontal and the
second term is vertical with respect to the splitting of $TE_b$.  Thus $(\Omega^B)^v 
 = R^E.$  Then
 for $k=1$ for simplicity, we have
\begin{eqnarray}\label{cw}
d_B\int_Z\tr(\Omega^B)^v  &=& \int_Z d_M\tr (R^E)
= \int_Z \tr \nabla^{\rm Hom}(R^E)\nonumber\\
&=&  \int_Z\tr [\nabla, R^E]\nonumber
= 0.
\end{eqnarray}
The equality $d_M \tr(R^E) =  \tr \nabla^{\rm Hom}(R^E)$ is an easy calculation using 
the fact that $R^E$ is skew-hermitian, and 
the last line uses the Bianchi identity.  

Thus we have constructed characteristic classes, and in particular a Chern character, 
 for the restricted class of
Bismut connections
without finding ${\rm Ad}_\G$-invariant functions on $\calg.$  Note that in the case of Remark
\ref{remark2}, these classes reduce to the classes discussed in Section 3, since for gauge transformations $\Omega = \Omega^v.$

The Chern character of the infinite rank superbundle associated to a family of coupled
Dirac operators is
\begin{equation}\label{FITrhs}
ch(\calE) = \int_Z ch(R^E) = \int_Z ch(\Omega^{S^+-S^-})\cup ch(\Omega^K).
\end{equation}
Since we have $ch(\Omega^{S^+-S^-})$ and not $\hat A(\Omega^M)$, this is not what we wanted!  

\medskip

Despite this failure, we now see what we can do with step (ii).
 We want to mimic the usual  ``cancellation of nonzero eigenspaces" in the heat equation proof of
 the index theorem, i.e., we want a connection that respects the splitting
 $$(S^+\otimes K)_b= {\ker}(\ddd^{\nabla^K}_b) \oplus \ker_+^\perp,$$
  where
$\ker_+^\perp$ equals $( {\ker}(\ddd^{\nabla^K}_b))^\perp$, and similarly for $S^-\otimes K.$  This perpendicular component is spanned by the eigensections of $\ddd^{\nabla^K} = \ddd_+$ with nonzero eigenvalues,
and $\ddd^{\nabla^K}$ is an isomorphism between these eigenspaces.  For this connection,
we expect that $ch(\ker_+^\perp) = ch(\ker_-^\perp)$ as forms 
computed with respect to this split connection.  This would imply
\begin{equation}\label{FITlhs} ch(\indd(\ddd^{\nabla^K})) = ch(\ker(\ddd^{\nabla^K}_+)) + 
ch(\ker_+^\perp) - ch( \ker(\ddd^{\nabla^K}_-))
 - ch(\ker_-^\perp)
= ch(\calE).
\end{equation}
This would finish (ii).  

The natural choice for such a connection is given by orthogonally projecting  the Bismut connection to 
the kernel and its perpendicular complement.  The problem is that the isomophism
$\ddd^{\nabla^K}_+: ({\rm ker}\ddd_+)^\perp\to ({\rm ker}\ddd_-)^\perp$
 is not a $\G$-isomorphism, since $\ddd ^{\nabla^K}_+$ is far from an element of $\G$.  
 However, we can replace  $D = \ddd^{\nabla^K}_+$ with its unitarization 
 $D^u = D/|D^*D|^{1/2}$ on these complements.  $D^u$ is a zero order invertible pseudodifferential operator, precisely the type of operator treated in Section 4.  
 
This motivates extending the structure group to a semidirect product 
$\tilde \G = \G\ltimes \pdos$, with $\G$ acting on $\pdos$ by conjugation.  We have to extend
our definition of characteristic classes from $\G$ to $\tilde \G$, but this is straightforward 
based on the earlier work:  for a $\tilde \G$-connection with curvature $\tilde \Omega = 
(V, s, B)\in \tilde {\calg}$ (so $B\in \pdoo$), we consider expressions like 
$$\int_{Z_b} \tr((s^{v^{\nabla_b}})^k)\dvol_{M/B} + \int_{Z_b}\tr ((\sigma_0(B))^k).$$
It is now straightforward to check that (\ref{FITlhs}) holds.  However, we are not claiming
that the Chern character form for this projected connection in (\ref{FITlhs}) is closed. Even 
if it is, we are definitely not claiming that it is cohomologous to the Chern character form in 
(\ref{FITlhs}), as this would give the wrong formula for the FIT.  

Despite the glaring problems with these arguments, we see
that when we try to prove the FIT directly on $\calE$, the extended structure group 
$\tilde \G$ naturally occurs.  

\begin{remark}  Recall that the starting point for Donaldson theory is the moduli space
$\calA/\calG$, where $\calA = \calA(F)$ is the space of connections on $F$ and $\calG = \calG(F)$ is the gauge group.  The action of $\calG$ on $\calA$, $g\cdot \nabla = 
g\nabla g^{-1}$, extends to an action of $\G$ by $(\phi, f)\cdot\nabla = f(\phi^{-1})^*\nabla
f^{-1}$, since for $\phi= Id,$ $f$ is a gauge transformation.

Thus it is natural to consider the moduli space 
 $\calA/ \G$.  However, this space seems to be a fat point in the following sense.
\medskip

\begin{conjecture} \cite{lprs} The orbit of a generic connection is dense.  
\end{conjecture}
\medskip 

An example of a nongeneric connection is a flat connection. As justification for the conjecture, 
it can be shown  that 
the normal space to the $ \G$-orbit $\calO_\nabla$ in the $L^2$ metric on 
$T_\nabla\calA = \Lambda^1(Z,\eend(F))$ is 
$$\{\alpha\in \Lambda^1(Z,\eend(F)): \nabla^*\alpha = 0\ {\rm and}\ R^F(\cdot, V)_z \perp \alpha_z, \forall
z\in Z, \forall V\in T_zZ\}.$$ (The equation $\nabla^*\alpha = 0$ is  the equation for 
the normal space to the gauge orbit of $\nabla$.) It is plausible  that for a generic connection, the curvature equation
$R^F(\cdot, V)_z \perp \alpha_z$ and its higher covariant derivatives 
form an overdetermined system of equations, and so has only the zero solution.  
We expect standard gauge theory techniques to help prove the conjecture.
\end{remark}

\bibliographystyle{amsplain}
\bibliography{Paper2}

\providecommand{\bysame}{\leavevmode\hbox to3em{\hrulefill}\thinspace}
\providecommand{\MR}{\relax\ifhmode\unskip\space\fi MR }
\providecommand{\MRhref}[2]{%
  \href{http://www.ams.org/mathscinet-getitem?mr=#1}{#2}
}
\providecommand{\href}[2]{#2}
\begin{thebibliography}{10}

\bibitem{adler}
{Adler, M.}, {van Moerbeke, P.}, and {Vanhaecke, P.}, \emph{{A}lgebraic
  {I}ntegrability, {P}ainlev\'e {G}eometry and {L}ie {A}lgebras}, Results in
  Mathematics and Related Areas. 3rd Series. A Series of Modern Surveys in
  Mathematics, {V}ol. 47, Springer-Verlag, Berlin, 2004.

\bibitem{atiyah}
{Atiyah, M.}, \emph{{\rm Circular symmetry and stationary phase
  approximation}}, {\it Ast\'erisque} \textbf{131} (1984), 43--59.

\bibitem{AS-IV}
{Atiyah, M.} and {Singer, I. M.}, \emph{{\rm The index of elliptic operators.
  {I}{V}.}}, {\it Annals of Math.} \textbf{93} (1971), 119--138.

\bibitem{BGV}
{Berline, N.}, {Getzler, E.}, and {Vergne, M.}, \emph{Heat {K}ernels and
  {D}irac {O}perators}, Grundlehren der Mathematischen Wissenschaften 298,
  Springer-Verlag, Berlin, 1992.

\bibitem{Bismut2}
{Bismut, J. M.}, \emph{\rm {The {A}tiyah-{S}inger index theorem for families of
  {D}irac operators: two heat equation proofs}}, {\it Inventiones Math.}
  \textbf{83} (1986), 91--151.

\bibitem{chern}
{Chern, S.-S.}, \emph{Complex {M}anifolds without {P}otential {T}heory},
  Springer-Verlag, New York, 1979.

\bibitem{C-S}
{Chern, S.-{S}.} and {Simons, J.}, \emph{{\rm Characteristic forms and
  Geometric Invariants}}, {\it Annals of Math.} \textbf{99} (1974), no.~1,
  48--69.

\bibitem{c-e}
{Chevalley, C.} and {Eilenberg, S.}, \emph{{\rm Cohomology groups of {L}ie
  groups and {L}ie algebras}}, {\it Trans. AMS} \textbf{63} (1948), 85--124.

\bibitem{dieu}
{Dieudonn\'e, J.}, \emph{A {H}istory of {A}lgebraic and {D}ifferential
  {T}opology 1900 -- 1960}, Birkh\"auser, Boston, 1989.

\bibitem{dupont}
{Dupont, J. L.}, \emph{Curvature and {C}haracteristic {C}lasses}, Lect. Notes
  Math. 640, Springer-Verlag, Berlin, 1978.

\bibitem{eells}
{Eells, J.}, \emph{{\rm A setting for global analysis}}, {\it Bull. Amer. Math.
  Soc.} \textbf{72} (1966), 751--807.

\bibitem{fgl}
{{Fedosov, B.}, {Golse, F.}, {Leichtnam, E.}, and {Schrohe, E.}}, \emph{{\rm
  The noncommutative residue for manifolds with boundary}}, {\it J. Funct.
  Analysis} \textbf{142} (1996), 1--31.

\bibitem{Freed}
{Freed, D.}, \emph{{\rm Geometry of Loop Groups}}, {\it J. Diff. Geom.}
  \textbf{28} (1988), 223--276.

\bibitem{g-m}
{{Garland, H.}, {Murray, M. K.}}, \emph{{\rm Kac-{M}oody algebras and periodic
  instantons}}, {\it Commun. Math. Phys.} \textbf{120} (1988), 335--351.

\bibitem{gdsw}
{{Gauntlett, J.P.}, {Martelli, D.}, {Sparks, J.}, {Waldram, D.}}, \emph{{\rm
  Sasaki-{E}instein metrics on ${S}^2\times {S}^3$}}, {\it Adv. Theor. Math.
  Phys.} \textbf{8} (2004), 711, hep--th/0403002.

\bibitem{guest}
{Guest, M.}, \emph{Harmonic {M}aps, {L}oop {G}roups, and {I}ntegrable
  {S}ystems}, LMS Student Texts, {V}ol. 38, Cambridge U. Press, Cambridge,
  1997.

\bibitem{Hamilton}
{Hamilton, R.}, \emph{{\rm Nash-{M}oser implicit function theorems}}, {\it
  Bull. Amer. Math. Soc.} \textbf{7} (1986), 65--222.

\bibitem{hatcher}
{Hatcher, A.}, \emph{Algebraic {T}opology}, Cambridge U. Press, Cambridge, UK,
  2002, www.math.cornell.edu/~hatcher/AT/ATpage.html.

\bibitem{huse}
{Husemoller, D.}, \emph{Fibre {B}undles, {\rm 1st ed.}}, Springer-Verlag, New
  York, 1966.

\bibitem{andres}
{Larrain-{H}ubach, A.}, \emph{{\rm Explicit computations of the symbols of
  order 0 and -1 of the curvature operator of ${\Omega} {G}$}}, {\it Letters in
  Math. Phys.} \textbf{89} (2009), 265--275.

\bibitem{lprs}
{Larrain-Hubach, A.}, {Paycha, S.}, {Rosenberg, S.}, and {Scott, S.}, in
  preparation.

\bibitem{lrst2}
{Larrain-Hubach, A.}, {Rosenberg, S.}, {Scott, S.}, and {Torres-Ardila, F.}, in
  preparation.

\bibitem{lrst}
\bysame, \emph{{\rm Characteristic classes and zeroth order pseudodifferential
  operators}}, {\it Spectral {T}heory and {G}eometric {A}nalysis}, Contemporary
  Mathematics, {V}ol. 532, AMS, 2011.

\bibitem{l-m}
{{Lawson, H. Blaine} and {Mickelsohn, M.}}, \emph{{\it Spin {G}eometry}},
  Princeton U. Press, Princeton, NJ, 1989.

\bibitem{L-N}
{Lesch, M.} and {Neira Jimenez, C.}, \emph{{\rm Classification of traces and
  hypertraces on spaces of classical pseudodifferential operators}},
  arXiv:1011.3238.

\bibitem{L-P}
{Lescure, J.-{M.}} and {Paycha, S.}, \emph{{\rm Uniqueness of multiplicative
  determinants on elliptic pseudodifferential operators}}, {\it Proc. London
  Math. Soc.} \textbf{94} (2007), 772--812.

\bibitem{MRT}
{Maeda, Y.}, {Rosenberg, S.}, and {Torres-Ardila, F.}, \emph{{\rm Riemannian
  geometry on loop spaces}}, arXiv:0705.1008.

\bibitem{milnor}
{Milnor, J.}, \emph{Characteristic {C}lasses}, Princeton U. Press, Princeton,
  1974.

\bibitem{mis-r-t}
{Misiolek, G.}, {Rosenberg, S.}, and {Torres-Ardila, F.}, in preparation.

\bibitem{m-v}
{{Murray, M. K.} and {Vozzo, R.}}, \emph{{\rm The caloron correspondence and
  higher string classes for loop groups}}, {\it J. Geom. Phys}. \textbf{60}
  (2010), 1235--1250.

\bibitem{Omori}
{{Omori, H.}}, \emph{Infinite-{D}imensional {L}ie {G}roups}, A.M.S.,
  Providence, RI, 1997.

\bibitem{paychasurvey}
{Paycha, S.}, \emph{{\rm Chern-{W}eil calculus extended to a class of infinite
  dimensional manifolds}}, arXiv:0706.2554.

\bibitem{ponge}
{Ponge, R.}, \emph{{\rm Traces on pseudodifferential operators and sums of
  commutators}}, ar{X}iv:0607.4265.

\bibitem{Seg}
{Pressley, A.} and {Segal, G.}, \emph{{L}oop {G}roups}, Oxford University
  Press, New York, NY, 1988.

\bibitem{sts}
{Reyman, A.G.} and {Semenov-Tian-Shansky, M.A.}, \emph{{I}ntegrable {S}ystems
  {II}: Group-{T}heoretical {M}ethods in the {T}heory of {F}inite-{D}imensional
  {I}ntegrable {S}ystems}, Dynamical systems. {VII}, Encyclopaedia of
  Mathematical Sciences, Vol. 16, Springer-{V}erlag, Berlin, 1994.

\bibitem{rochon}
{Rochon, F.}, \emph{{\rm Sur la topologie de l'espace des op\'erateurs
  pseudodiff\'erentiels inversibles d'ordre 0}}, {\it Ann. Inst. Fourier}
  \textbf{58} (2008), 29--62.

\bibitem{rosenberg}
{Rosenberg, S.}, \emph{The {L}aplacian on a {R}iemannian {M}anifold}, Cambridge
  U. Press, Cambridge, UK, 1997.

\bibitem{scott}
{Scott, S.}, \emph{Traces and {D}eterminants of {P}seudodifferential
  {O}perators}, Oxford U. Press, Oxford, 2010.

\bibitem{witten}
{Witten, E.}, \emph{{\rm Quantum field theory and the {J}ones polynomial}},
  {\it Commun. Math. Phys.} \textbf{121} (1989), 351--399.

\end{thebibliography}

\end{document}